\documentstyle[12pt]{article}
\title{ On The Eigenvalues of Some Vectorial \\
        Sturm-Liouville Eigenvalue Problems \\
        }
\author{ Hua-Huai Chern }
\date{ March 25, 1998 }
\begin{document}
\baselineskip=20pt
\maketitle
\vskip -0.5cm
\centerline{\footnotesize Department of Mathmatics}
\centerline{\footnotesize National Chung Cheng University}
\centerline{\footnotesize Minghsiung, Chiayi 621 }
\centerline{\footnotesize Taiwan }
\centerline{\footnotesize E-mail: hhchern@math.ccu.edu.tw }
\vskip 0.5cm
\par
\vskip 2.5cm
\abstract
The author tries to derive the asymptotic expression of the 
large eigevalues of some vectorial Sturm-Liouville 
differential equations. A precise description for the 
formula of the square root of the large eiegnvalues up to  
the $O(1/n)$-term is obtained.
  \par
\vskip 0.5cm
\noindent {\bf AMS(MOS) subject classifications. } 34A30, 34B25. \par
\vskip 1cm
\noindent Abbreviated title: Asymptotics of eigenvalues of 
vectorial S-L problems
 \newpage
%
%
%

\section{Introduction}
For the scalar case, it is well-known that: if $q^{(m)} (x) \in {\cal L}^1 
(0, \pi ) $, the equation
\[
   -y'' +q(x) y = \lambda y, \quad y' (0)-hy(0)=0, \quad y' (\pi )+H y (\pi ) 
=0 \] 
has eigenvalues possessing the following asymptotic form
\[
  \sqrt \lambda_n =n + \frac {a_\circ }n + \cdots +
         \frac {a_{[m/2]}}{n^{2[m/2]+1}} + \frac {\gamma_n}{n^{2[m/2]+1}},
\]
where $ a_\circ = (h+H+1/2 \int_0^\pi q(t)dt )/\pi $ and the other $a_j$, 
$ 1 \le j \le [ m /2 ] $ are constant depending on $h, H, q(x) $ and the 
higher order derivatives of $q(x)$ up to order $m$. The main technique for
deriving above expression is by using iterative method. One can find the
standard approaches in the appendix in \cite{LG} or in \cite{Hoch}. For the
vectorial case, the author only have seen analogous results in \cite{CA} for
the case $d=2$ and the boundary condtions are different to the cases the author
has considered. Motivated by the arguments in \cite{CA},
we still can use iterative method analogous to the one in
\cite{LG}, but only partial result as {\bf Theorem \ref{Th1}} can be obtained.
The technique we use in this section is similar to the one in \cite{CA}, but
the major difference between theirs and the author's is that the author locate
the eigenvalues by using the associated matrix differential equation and by
the theory of complex operator-valued functions. Hence, we can dispose any
dimensional case.
Consider the following vectorial Sturm-Liouville differential systems
\begin{equation} \label{VSL}
  -\phi '' (x) + P(x)\phi (x)= \lambda \phi (x),  \quad
   \phi'(0)+H_L \phi (0) = {\bf 0 } = \phi ' (\pi )+H_R \phi (\pi)
\end{equation}
where $P(x)$ is an $N\times N$ real symmetric matrix-valued function and $H_L , H_R $ are $N\times N $ real symmetric constant matrices

Consider the matrix differential equations related to (\ref{VSL}) and (FS),
respectively,
\begin{eqnarray*}
    -Y'' (x) +P(x) Y (x)= \lambda Y(x), \quad Y(0)=I, \quad Y'(0)=-H_L  ,\\
    -{\cal C}'' (x) = \lambda {\cal C} (x), \quad {\cal C}(0)=I, \quad {\cal C}'(0)=-H_L,
\end{eqnarray*}
the solutions $Y(x; \lambda )$ and $ {\cal C} (x;\lambda ) $ are connected by
the following identity
\begin{equation} \label{1.2}
  Y(x; \lambda )={\cal C} (x;\lambda)+ \int_0^x {\cal K}( x,t){\cal C} (t;\lambda ) dt, \quad \forall \lambda
\end{equation}
and ${\cal C} (x;\lambda)= \cos (\sqrt{\lambda}x)I-\sin (\sqrt{\lambda }x)H_L /\sqrt{\lambda}  $. The connection (\ref{1.2}) is easy to prove, one
can see a more general version in \cite{CS1} .

In this paper, we also use the symbol $(P,H_L ,H_R )$ to denote the
eigenvalue problem (\ref{VSL}). 
Denote $\Sigma (P,H_L ,H_R ) $ the sequence of
eigenvalues of $(P,H_L ,H_R )$. We can arrange those elements in it in
ascending order as
\[
      \mu_0 \le \mu_1 \le \ldots \le \mu_k \le \ldots .
\]
We use the symbol $\sigma (P, H_L, H_R )$ to represent the set of 
eigenvalues of $(P,H_L, H_R )$, and arrange those elements in this set
in ascending order as
\[
      \lambda_0 < \lambda_1 < \ldots < \lambda_k < \ldots .
\]
The multiplicity of each eigenvalues $\lambda_k$ is denoted by $m_k$. For
each positive integer $n$, an index set $\Lambda_n$ is defined as
\[
         \Lambda_n = \{ k: \lambda_k^{ 1/2} \in B_{ 1/4} (n) \}
\]
where $B_r(\xi ) = \{ z: |z-\xi| < r \}$.

With above basic terminology and notations, we shall give the proof of
the following theorem in section 2:
\newtheorem{1}{Theorem}\begin{1} \label{Th1} Suppose $\lambda_k \in \sigma 
(P,H_L, H_R )$ with multiplicity $m_k$. For $n$ is sufficiently large,
if $k \in \Lambda_n $, then
\begin{equation} \label{1.3}
  \lambda_k^{\frac 12}= n + \frac {a_k}n + \gamma_n (k),
\end{equation}
where $ a_k $ is the characteristic value of the $N \times N $ real symmetric
matrix $ \frac 1\pi ( \frac 12 \int_0^\pi P(t)dt + H_R-H_L )$, and 
$\gamma_n (k)=o(1/n^2)$. Furthermore,
$ \sum_{k \in \Lambda_n} m_k = N$, where $N$ is the dimension of the system
(\ref{VSL}).
\end{1}

%
%

\section{Asymptotic analysis of eigenvalues}

In order to make this section become more readable, the author outlines the
sketch of the whole approach as below: The eigenvalues $\lambda_k $, $k \ge 0$
of (\ref{VSL}) can be locate by determining whether the matrix-valued function
\[
   W(\lambda_k ) = Y' (\pi ; \lambda_k ) + H_R Y (\pi ; \lambda_k )
\]
is singular or not. We use iterative method so as to write
\begin{equation} \label{2.0}
    W(\lambda ) = \Psi (\lambda ) + {\cal E} (\lambda )
\end{equation}
where we shall see that $\Psi (\lambda )$ has a quite neat and simple form from
which we can determine those values $\tilde{\lambda }$ making $\Psi (\lambda )$ 
be singular. By the extension theorem of Rouch\'e's theorem on operator-valued
functions, we are getting close to locate the eigenvalues of (\ref{VSL}), 
although we can not locate them completely and explicitly as we did in the  
scalar case. With above short sketch in mind, we now start the lengthy and
tedious computations on the asymptotics of eigenvalues of (\ref{VSL}).

Rewrite
\begin{equation} \label{2.1}
  Y( x; \lambda ) = {\cal C} ( x ; \lambda ) +
 \int_0^x \frac {\sin \sqrt{\lambda } (x-t) }{\sqrt{\lambda }} P(t) 
 Y(t; \lambda ) dt. \end{equation}
Then
\begin{eqnarray}  \nonumber
     Y_1 (x; \lambda ) &=&  Y (x; \lambda ) - {\cal C} ( x ; \lambda ) \\
                       &=&
 \int_0^x \frac {\sin \sqrt{\lambda } (x-t) }{\sqrt{\lambda }} P(t) 
 {\cal C} (t; \lambda ) dt  +
 \int_0^x \frac {\sin \sqrt{\lambda } (x-t) }{\sqrt{\lambda }} P(t) 
 Y_1 (t; \lambda ) dt \label{2.2} \\
    &=& I_1 +  \nonumber
 \int_0^x \frac {\sin \sqrt{\lambda } (x-t) }{\sqrt{\lambda }} P(t) 
 Y_1 (t; \lambda ) dt  .
\end{eqnarray}
Using integration by part, we can compute $I_1 $ as
\[
   I_1 = \frac {\sin \sqrt{\lambda }x }{\sqrt{\lambda }} {\cal K}(x,x) 
        + \frac {\cos \sqrt{\lambda }x}{\lambda }( {\cal K} (x,x) H_L 
          + \frac 14 (P(x)- P(0))) + O (|\lambda |^{- \frac 32 } ).
\]
Combining with (\ref{2.2}), we have 
\begin{eqnarray*}
 Y_1 (x; \lambda ) &=& \frac {\sin \sqrt{\lambda }x }{\sqrt{\lambda }} {\cal 
                       K}(x,x) +
                      \frac {\cos \sqrt{\lambda }x}{\lambda }( {\cal K} (x,x)
                       H_L + \frac 14 (P(x)- P(0))   \\
                   & & \mbox{} +
 \int_0^x \frac {\sin \sqrt{\lambda } (x-t) }{\sqrt{\lambda }} P(t) 
 Y_1 (t; \lambda ) dt  + O (|\lambda |^{- \frac 32 } ).
\end{eqnarray*}
As $ |\Im \lambda | < \kappa $ for some constant $\kappa $ fixed, by Gronwall's
inequality, we can prove that $ || Y_1 (x; \lambda ) ||_\infty = O( 
|\lambda |^{-1/2} )$, as $ |\lambda | $ is sufficiently large. From above 
argument, we have 
\begin{eqnarray} \label{2.3} 
 Y( x ; \lambda ) &=& \cos \sqrt{\lambda } x 
                     ( I + \frac 1{\lambda } {\cal K} (x,x) H_L + \frac 1{4 
\lambda }(P(x)-P(0))  \\
& & \mbox{}+ \frac {\sin \sqrt{\lambda } x}{\sqrt{\lambda }}
({\cal K}(x,x)-H_L ) + O (|\lambda |^{- \frac 32 } ). \nonumber
\end{eqnarray}
On the other hand, from (\ref{2.1}),
\begin{equation} \label{2.4}
  Y '( x; \lambda ) = {\cal C}' ( x ; \lambda ) +
 \int_0^x \cos \sqrt{\lambda } (x-t) P(t) 
 Y(t; \lambda ) dt. \end{equation}
Plugging (\ref{2.3}) into (\ref{2.4}), and repeating the analysis as above,
we have
\begin{eqnarray} \label{2.5}
  Y' (x; \lambda )&=& - \sqrt{\lambda } \sin \sqrt{\lambda }x I
  + \cos \sqrt{\lambda }x ( {\cal K} (x,x) -H_L )  \nonumber \\
  & &\mbox{} + \frac {\sin \sqrt{\lambda }x}{\sqrt{\lambda }}
   (\frac {P(x)+P(0)}2 + \frac 12 \int_0^x P(t){\cal K}(t,t ) dt
     -{\cal K} (x,x) H_L )  \\
  & &\mbox{} + \frac {\cos \sqrt{\lambda }x }{\lambda }
      ( \frac 14 (P' (x) - P'(0) ) + \frac 12 ( \int_0^x P(t) {\cal K}
        (t,t) dt )H_L  \nonumber \\
  & & \mbox{} +\frac 18 \int_0^x (P^2 (t) -P(t)P(0))dt
              +\frac 12 (P(x){\cal K} (x,x) -H_L ) )+ O(|\lambda |^{-\frac 32 
})   \nonumber
\end{eqnarray}
By (\ref{2.3}) and (\ref{2.5}), we have
\begin{equation} \label{2.6}
   W( \lambda )= -\sqrt{\lambda } \sin \sqrt{\lambda } \pi I + \cos 
\sqrt{\lambda} \pi {\cal G}_1 + \frac {\sin \sqrt{\lambda } \pi }{ 
\sqrt{\lambda}} {\cal G}_2 +\frac {\cos \sqrt{\lambda } \pi }{\lambda } 
{\cal G}_3 +O(|\lambda |^{-\frac 32}),
\end{equation}
where
\[
   {\cal G}_1 = H_R -H_L + {\cal K} (\pi , \pi ) , \]
\begin{eqnarray*}
   {\cal G}_2 &=& \frac {P(\pi )+ P (0)}2 + \frac 12 \int_0^\pi P(t){\cal K}   
   (t,t) dt + H_R {\cal K} (\pi ,\pi )\\
     & &\mbox{ }-{\cal K} (\pi , \pi ) H_L -H_R H_L ,
\end{eqnarray*}
and
\begin{eqnarray*}
  {\cal  G}_3 &= &\frac {P' (\pi ) - P'(0)}4+ \frac 12 ( \int_0^\pi P(t) {\cal 
K} (t,t) dt )H_L + \frac 14 H_R (P(\pi ) -P(0)) \\
              & & \mbox{} + \frac 18 \int_0^\pi (P^2 (t) -P(t)P(0))dt
                +\frac 12 (P(\pi ) {\cal K} (\pi , \pi ) -H_L )+
                 H_R {\cal K}(\pi , \pi ) H_L
\end{eqnarray*}
Write $\lambda = \mu^2 $. Then the term $ \Psi (\lambda ) $ in  (\ref{2.0}) 
is defined as $\Psi (\mu^2 ) = -\mu \sin \mu \pi I + \cos \mu \pi {\cal G}_1 $.
The remainder terms of (\ref{2.6}) is denoted by ${\cal E} (\mu ^2 ) $. Denote
$\varpi (\mu )=\det ( W( \mu^2 ))$ and $ \mu = \sigma + it $. By the following 
identities 
\begin{eqnarray*}
  \sin (\sigma \pi + i t \pi ) &=& \sin \sigma \pi \cosh t \pi
    + i \sinh t\pi \cos \sigma \pi , \\
  \cos (\sigma \pi + i t \pi ) &=& \cos \sigma \pi \cosh t \pi
    - i \sinh t\pi \sin \sigma \pi
\end{eqnarray*}
and using the Laplace expansion of determinants, we have the following result:
for any complex $\mu $, $ \varpi( \mu ) = \vartheta ( \mu )+ \varepsilon ( \mu ) $
where $ \vartheta (\mu ) = \mu ^N \sin ^N \mu \pi $ and $ | \varepsilon (\mu 
) / \vartheta (\mu ) | = O ( 1/|\mu | ) $. With above facts and Rouch\'e's 
theorem, there lie exactly $N$ zeroes of $\varpi (\mu ) $ 
within a suitable neighborhood of any sufficiently large integer $n$. If we
denote a zero of $\varpi (\mu ) $ by $ n+ \delta_{ni} $, it is not difficult to 
see that $\delta_{ni}=O(1/n) $. In fact, we can compute 
$ \lim_{n \rightarrow \infty }\delta_{ni} $, and hence we have the first part
of {\bf Theorem \ref{Th1} }.

To touch down the goal, let us back to (\ref{2.0}). The matrix-valued function 
$\Psi (\mu )$ actually provides us much information. Since the
matrix ${\cal G}_1$ is real symmetric, there exists an orthogonal matrix $U$,
$U^{-1}=U^*$, such that
\[
   U^* \Psi (\mu ) U = \mbox{diag}[ \varrho_1 (\mu ), \varrho_2 (\mu ),
 \cdots, \varrho_N (\mu )] 
\]
where $ \varrho_i (\mu )= -\mu \sin \mu \pi +\sigma_i \pi \cos \mu \pi $, $ 1 
\le i \le N$, and $\sigma_i , 1 \le i \le N $ are those characteristic values 
of constant matrix ${\cal G}_1 $. By (\ref{2.0}) and (\ref{2.6}), we have
\begin{eqnarray} \label{2.7}
   U^* W (\mu^2 ) U &=& U^* \Psi (\mu^2 ) U+ U^* {\cal E} (\mu^2 ) U \\
                    &=& U^* \Psi (\mu^2 ) U + \frac {\sin \mu \pi}{\mu }
                        U^* {\cal G}_2 U + O (\frac 1{\mu^2 }). \nonumber
\end{eqnarray}
Using the following Rouch\'e's Theorem for operator-valued functions in 
\cite{GGK}, we can suitably locate the eigevalues of (\ref{VSL}) from (\ref{2.7}). 

\newtheorem{2.1}[1]{Theorem}\begin{2.1} \label{Th2.1}                      
Let $W(\mu )$ and $S(\mu )$ be ${\cal L 
(H)}$-valued holomorphic function defined on a region $\Omega $ enclosed by a
contour $\Gamma $. Suppose $W$ is normal with respect to $\Gamma $ and let $V=
W+S$. If
\begin{equation} \label{2.8}
           || W^{-1} (\mu ) S (\mu ) || < 1 ,
\end{equation}
then $V$ is also normal with respect to $\Gamma $, and 
\[
    m (\Gamma ; V(\cdot ))= m (\Gamma ; W(\cdot )).
\]
\end{2.1}

One can refer to \cite[ chapter XI, \S 9, 205--211]{GGK} for the definition 
of all the notations in above 
theorem. In our case, the first task we must face is the estimate of the norm 
of $ U^* \Psi^{-1} (\mu ) {\cal E} (\mu )U $ on some suitable contours. The 
estimate will be taken with respect to $\infty$-norm, but this task is very 
tedious, we put the details in the appendix and give only the statement as 
below here. In the following lemma, the contour is chosen as: for any complex
number $\mu_\circ $ ,
\begin{equation} \label{2.9}
          \Gamma_n (\mu_\circ ) = \{ \mu : |\mu - \mu_\circ | = O (\frac 
1{n^2}). \} 
\end{equation}
in which we do not specify the radius since we shall obtain a uniform estimate.

\newtheorem{2.2}[1]{Lemma}\begin{2.2} \label{Th2.2}
Let $\mu_\circ $ be a zero of the transcendental equation
\[
    \mu \sin \mu \pi - \alpha \pi \cos \mu \pi =0,
\]
locating in a suitable neighborhood of some sufficiently large integer 
$n$, where $\alpha \pi $ is a characteristic value of ${\cal G}_1 $. Let $ 
\Gamma_n (\mu_\circ ) $ is chosen as (\ref{2.9}) with the same $n$. Then 
\[
  || U^* \Psi^{-1} (\mu ) {\cal E} (\mu )U ||_\infty =O(\frac 1n), \quad 
\forall \mu \in \Gamma_n (\mu_\circ ) . \]
\end{2.2}
Using Largrange inversion formula, it is not 
difficult to see that the zero of the transcendal equation
\[
    \mu \sin \mu \pi - \alpha \pi \cos \mu \pi =0,
\]
lying in a suitable neighborhood of positive integer $n$ can be expressed 
as 
\[
          \mu_n = n+ \frac {\alpha }n + \frac {\kappa (\alpha )}{n^3 } + 
\cdots . 
\]
By {\bf Theorem \ref{Th2.1}} and {\bf Lemma \ref{Th2.2}}, if $ \alpha \pi $ is
a characteristic value of ${\cal G}_1 $, then there exist at least one 
eigenvalue $ \lambda_n ( \alpha ) $ of (\ref{VSL}) whose square root $\sqrt
\lambda_n (\alpha ) $ satisfies
\begin{equation} \label{2.10}
      \sqrt{\lambda_n} (\alpha ) = n + \frac {\alpha }n+ \gamma_n (\alpha ),
\end{equation}
where $\gamma_n (\alpha )=o(\frac 1{n^2} )$. We summarize above as
\newtheorem{2.3}[1]{Lemma}\begin{2.3} \label{Th2.3}
Suppose the eigenvalues of (\ref{VSL}) satisfying (\ref{2.10}) in asymptotic 
sense. Then 
\[
\{ \alpha \pi : \mbox{all possible $ \alpha $  that (\ref{2.10}) 
holds } \} \supseteq \sigma ({\cal G}_1 ), 
\]
 where $ \sigma ({\cal G}_1 ) $ is 
the spectral set of ${\cal G}_1 $.
\end{2.3}
 On the other hand, we have the following 
result. 

\newtheorem{2.4}[1]{Lemma}\begin{2.4} \label{Th2.4} If any sequence of eigenvalue of
(\ref{VSL}) has the asymptotic expression as (\ref{2.10}), then $\alpha \pi 
\in \sigma ({\cal G}_1 ) $ and $ \gamma_n (\alpha )= o(1/n^2) $.
\end{2.4}

\noindent {\bf Proof. } By (\ref{2.6}), 
\begin{equation} \label{2.11}
  W(\mu^2 ) = -\mu \sin \mu \pi I + \cos \mu \pi {\cal G}_1 +\frac {\sin \mu 
\pi}{\mu } {\cal G}_2 + \frac {\cos \mu \pi }{\mu^2 } {\cal G}_3 
+ O (\frac 1{|\mu |^3 }).
\end{equation}
Let $ \mu_n = \sqrt \lambda_n (\alpha ) $ as given by (\ref{2.10} ), then 
plugging $\mu_n $ into (\ref{2.11}), we have 
\begin{eqnarray}  \label{2.12}
   W( \mu_n^2 ) &=& - (n+\frac {\alpha }n + \gamma_n (\alpha ) )
   \sin (n+ \frac {\alpha }n + \gamma_n (\alpha ) )\pi I  \nonumber
   +\cos (n+ \frac {\alpha }n + \gamma_n (\alpha ) )\pi {\cal G}_1 \\
   & &\mbox{ }+ \frac 1{n+ \frac {\alpha }n + \gamma_n (\alpha ) }
   \sin (n+ \frac {\alpha }n + \gamma_n (\alpha ) )\pi {\cal G}_2   \\
   & &\mbox{}+\frac 1{(n+ \frac {\alpha }n + \gamma_n (\alpha ))^2 }\cos (n+ 
\frac {\alpha }n \gamma_n (\alpha ) )\pi {\cal G}_3+ O(\frac 1{n^3}). \nonumber
\end{eqnarray}
Since 
\begin{eqnarray*}
  \sin (\mu _n \pi )  
      &= & (-1)^n (\frac {\alpha }n + \gamma _n (\alpha ) )\pi + O(\frac 1{n^3 
}), \\   \cos(\mu _n \pi )
      & = & (-1)^n ( 1- \frac 1{2} (\frac {\alpha }n + \gamma _n (\alpha ) )^2 
\pi ^2 ) + O(\frac 1{n^4 }).
\end{eqnarray*}
Substituting them into (\ref{2.12}) and equating the terms by comparing orders,
we have
\begin{equation} \label{2.13}
 W( \mu _n^2 ) = (-1)^n [ \alpha \pi I- {\cal G}_1 +n\gamma _n (\alpha ) \pi I
 ]+ O (\frac 1{n^2} ).
\end{equation}
Since $W( \mu^2_n ) $ is singular for all $n$ sufficiently large, and $n \gamma 
_n (\alpha ) \rightarrow 0 $, by the continuity of determinant function, the 
matrix $  \alpha \pi I-{\cal G}_1 $ must be
singular. Otherwise, letting $ {\cal T} = \alpha \pi I - {\cal G}_1 $, we have,
by (\ref{2.13}), that 
\[
    || W(\mu_n^2 ) - (-1)^n {\cal T} || \le |n \gamma_n (\alpha ) \pi |+  
     O(\frac 1{n^2}).
\]
Thus if ${\cal T} $ is invertible, then so is $W(\mu_n^2 )$ for $n$ 
sufficiently large, which is absurd.  Therefore, 
 $\alpha \pi $ is a characteristic value of ${\cal G}_1
$.

By {\bf Lemma \ref{Th2.3} } and {\bf Lemma \ref{Th2.4}}, we have {\bf Theorem 
\ref{Th1} }. 

\noindent {\bf Remark}. Above approaches can be extended to some other type of
boundary conditions, the author has disposed some kinds of them in 
\cite{Chern}, e.g., the Dirichlet boundary conditions. The determination on the 
first term in (\ref{1.3}) is quite precise, while 
on the second or higher terms, there are still some ambiguities need to be 
excluded.

%
%

\par
\vskip 1cm
\noindent{\Large \bf Appendix. 
The estimate of $|| U^* \Psi ^{-1} (\mu ) {\cal E} (\mu 
) U ||_\infty $  }

\noindent In section 2, we give the statement on the estimate of 
$|| U^* \Psi ^{-1} (\mu ) {\cal E} (\mu ) U ||_\infty $. The proof of
that lemma lies on the estimate of orders on the prescribed contour.

By (\ref{2.6}) and (\ref{2.7}),
\begin{eqnarray*}    \label{a.1}
 U^* \Psi ^{-1} (\mu ) {\cal E} (\mu ) U &=&
  {\cal R} (\mu ) ( \frac {\sin \mu \pi }{\mu } U^* {\cal G}_2 U
     + \frac {\cos \mu \pi }{\mu^2 } U^* {\cal G}_3 U ) \\
   & &\mbox{}+ {\cal R} (\mu ) O(\frac 1{ |\mu |^3 }) ,
\end{eqnarray*}
where ${\cal R} (\mu ) =$ diag$[ \varrho_1 (\mu ), \ldots , \varrho_N (\mu ) 
] $. Let $\mu_\circ $ be a value making $\Psi (\mu_\circ )$ be singular. Then 
$\mu_\circ $ is a zero of $\varrho _j (\mu ) $ for some $ j$. According to the
distribution of the zeroes of the transcendental equation $\varrho_j (\mu )$,
we may assume $\mu_\circ $ lie in a suitable neighborhood of a sufficient large
integer $n$, and can be espressed as 
\[
   \mu_\circ = n + \frac an+ \frac b{n^3} + O(\frac 1{n^5}),
\]
where $b$ depends on $a$. With this kind $\mu_\circ $ as center, we choose
$\Gamma_n (\mu_\circ ) $ as defined in (\ref{2.9}). We shall show that: On
such contours, the two quantities
\[
    || {\cal R} (\mu ) \frac {\sin \mu \pi }{\mu } ||_\infty , \quad 
    || {\cal R} (\mu ) \frac {\cos \mu \pi }{\mu^2 } ||_\infty
\]
are of order $O(1/n )$, and hence the remainder term is of order less then
$O(1/n)$.

Here, we only give the demostration on the first quantity. As for the other
one, one can use the arguements in the coming paragraphs to prove it. 

Write 
\[
    {\cal R}(\mu) \frac {\sin \mu \pi }{\mu } =
   \left ( \begin{array}{cccc}
      \chi_1 (\mu ) &0 &\cdots &0 \\
       0 &\chi_2 (\mu ) & \cdots &\vdots \\
       \vdots & \vdots & \ddots &0 \\
         0  & \cdots & 0 &\chi_N (\mu )
    \end{array} \right )
\]
where $\chi_j (\mu ) = \sin \mu \pi /\mu (\mu \sin \mu \pi - \sigma_j \pi \cos 
\mu \pi ) $, $ 1 \le j \le N$. Then
\[
  || {\cal R} (\mu ) \frac {\sin \mu \pi }{\mu } ||_\infty
          = \max_{1 \le j \le N} || \chi_j (\mu ) ||_\infty
          = \max_{1 \le j \le N} \sup_{\mu \in \Gamma_n (\mu_\circ ) }
                             | \chi_j (\mu )|.
\]
We only to give the estimate on $|\chi_j (\mu ) |$ in the sense of growth 
of order. Write $\mu = s+ i t $, then for any $\mu \in \Gamma_n (\mu_\circ ) $, 
we have
\[
  \mu = n+ \frac an + \frac {\delta \cos \phi  }{n^2 }
               + i \frac {\delta \sin \phi }{n^2 } +\frac b{n^3 }+
O(\frac 1{n^5 }), \]
where $\delta $ is the radius of the contour and $0\le \phi \le 2\pi $. Then,
\begin{eqnarray}   \label{a.2}
   s &= & n+ \frac an + \frac {  \delta \cos \phi }
       {n^2 }+ \frac b{m^3 }+O(\frac 1{n^5 }), \\
    t  &= & \frac {\delta \sin \phi  }{n^2 }. \label{a.3}
\end{eqnarray}
The following identities are needed:
\begin{eqnarray*}
    |\sin \mu \pi |^2 &=& \frac 12 ( \cosh 2t \pi - \cos 2 s \pi ),
\\
   |\cos \mu \pi |^2 &=& \frac 12 ( \cosh 2t \pi + \cos 2 s \pi ),
\end{eqnarray*}
\[ 2\Re \{ \mu \sin \mu \pi \cos \bar{\mu } \pi  \}
     = s \sin 2 s \pi  + t \sinh 2 t\pi , \]
and
\begin{eqnarray*} \lefteqn{
  |\mu \sin \mu \pi  - \theta \pi \cos \mu \pi  |^2
      =|\mu \sin \mu \pi |^2 }  \\
      & &\mbox{ }+ (\theta \pi )^2 |\cos \mu \pi |^2
           -2\theta \pi \Re \{ \mu \sin \mu \pi \cos \bar{\mu } \pi  \},
\end{eqnarray*}
where $\theta \in {\bf R}$. Plugging (\ref{a.2}), (\ref{a.3}) into above 
identities, we have
\begin{eqnarray*}
|\sin \mu \pi |^2 & = &
   \frac {a ^2 \pi ^2}{n^2 } +\frac {2\pi^2}{n^3 }a \delta
\cos \phi    \\ & &
\mbox{}+\frac {\pi ^2}{3n^4 }(  -a^4 \pi ^2 + 6ab +3\delta ^2 )
    )+O(\frac 1{n^5 }) ,\quad \mbox{ if $a \neq 0$ ,}
\end{eqnarray*}
\[
|\sin \mu \pi |^2 =
 \frac {\pi ^2 \delta ^2 }{n^4 } + \frac {\pi ^4 }{3n^8 } \delta ^4 (1-2\cos
^2 \phi ) + O(\frac 1{n^9 }), \quad \mbox{ if $ a =0 $.}
\]
\begin{eqnarray*}
|\cos \mu \pi |^2 &= &
  1- \frac {a ^2 \pi ^2}{n^2 } -\frac {2\pi ^2}{n^3 } a \delta
\cos \phi \\ & &\mbox{}+\frac {\pi ^2}{3n^4 }
    (3\delta ^2
     (  a ^4 \pi ^2 - 6ab +1-2 \cos ^2 \phi  )
    )+O(\frac 1{n^6 }) , \quad \mbox{if $a \neq 0$,}
\end{eqnarray*}
\[ |\cos \mu \pi |^2 =
  1+\frac {\pi ^2}{n^4 } ( 1-2\cos ^2 \phi )
+ O(\frac 1{n^8 }), \quad \mbox{ if $ a =0 $.}
\]
Hence we have :
\begin{eqnarray*} \lefteqn{
   |\mu \sin \mu \pi  |^2  + (\theta \pi )^2 |\cos \mu \pi |^2 
      = 
    (\theta ^2 +a ^2)\pi ^2  +  \frac {2\pi ^2}n
    a \delta \cos \phi } \\ & &\mbox{ }+ \frac {\pi ^2 }{3n^2 } (3\delta ^2
     -a ^4 \pi ^2 + 6a ^3 -3(\pi a \theta )^2+6ab )
        +O(\frac 1{n^3 }), \quad  \mbox{ if $a \neq 0 $, }
\end{eqnarray*}
\[
 |\mu \sin \mu \pi  |^2  + (\theta \pi )^2 |\cos \mu \pi |^2 =
    (\theta \pi )^2 + \frac {\pi ^2 \delta ^2 }{n^2 }
                     + O(\frac 1{n^4 }), \quad
                     \mbox{ if $a = 0 $. }
 \]
On the other hand, we have
\begin{eqnarray*}
  2\Re \{ \mu \sin \mu \pi  \cos \bar{ \mu } \pi  \}
 &=&
  2a \pi + \frac {2\pi }{n } \delta \cos \phi  \\ & &\mbox{}
+\frac {2\pi }{3n^2} (3a ^2 +3b -2\pi ^2 a ^3 ) +
O(\frac 1{n^3 }) , \quad \mbox{ if $ a \neq 0$, }
\end{eqnarray*}
\[
  2\Re \{ \mu \sin \mu \pi  \cos \bar{ \mu } \pi  \} =
  \frac {2\pi }{n } \delta \cos \phi  + O(\frac 1{n^4 }),
      \quad \mbox{ if $a=0 $. }
\]
Therefore, we obtained
\begin{eqnarray*} \lefteqn{
  |\mu \sin \mu \pi  -\theta \pi \cos \mu \pi |^2 }\\
  &= &
   (a \pi -\theta \pi )^2  + \frac {2\pi ^2}n \delta (a -\theta)
   \cos \phi  +
   \frac {\pi ^2}{3n^2 } ( 6(a -\theta )(a ^2 + b)\\
    & &\mbox{}+3\delta ^2 -\pi^2 a ^2 (a ^2
                         +3\theta ^2 -4a \theta ))
     + O(\frac 1{n^ 3} ),
  \quad \mbox{ if $a \neq 0 $,}
\end{eqnarray*}
\[
|\mu \sin \mu \pi  -\theta \pi \cos \mu \pi |^2  =
  (\theta \pi )^2 - \frac {2\pi^2 }n \theta \cos \phi +\frac {\pi ^2 \delta
^2 }{n^2 }+O(\frac 1{n^4 }), \quad \mbox{ if $ a =0 $. }
\]
Finally, we have : On $\Gamma _n (\mu _\circ )$,
\begin{eqnarray*}
  |\sin \mu \pi  | &= & \mbox{}
   \frac 1n |a| \pi +o(\frac 1n ), \quad \mbox{ if
      $a \neq 0 $,} \\ 
  |\sin \mu \pi  | &= & \mbox{}
    \frac 1{n^2 } \delta \pi +o(\frac 1{n^2 }), \quad \mbox{ if $a =0$.}
\end{eqnarray*}
and
\begin{eqnarray*}
  |\mu \sin \mu \pi  -\theta \pi \cos \mu \pi | &=& \mbox{}
     O(1), \quad \mbox{ if $ \theta \neq a $.} \\
  |\mu \sin \mu \pi  -\theta \pi \cos \mu \pi | &=& \mbox{}
      \frac 1n \delta \pi +o(\frac 1n ), \mbox{if $ a \neq 0$, }
       \mbox{ if $ \theta = a $. }  \\
  |\mu \sin \mu \pi  -\theta \pi \cos \mu \pi | &=& \mbox{}
     \frac 1n \delta \pi +o(\frac 1n ),  \mbox{ if $a=0 $, }
     \mbox{ if $ \theta = a $. }
\end{eqnarray*}
So, if $ a \ne 0$, we obtain
\begin{eqnarray*}
 |\chi_j (\mu ) | &= &
          \frac {\frac 1n |a|\pi   + o(\frac 1n ) }{
                (n+o(1))(|\theta - a | +o(1)) } = O(\frac 1{n^2 }),
           \quad  \mbox{ if $ \theta \ne a $. }  \\
 |\chi_j (\mu ) | &= &
          \frac {\frac 1n |a|\pi   + o(\frac 1n ) }{
                (n+o(1))(\frac 1n \pi \delta +o(\frac 1n ))}
                 = \frac 1{n \delta }
                |\alpha \pi |+o(\frac 1n ),  \mbox{ if $ \theta = a $. }
\end{eqnarray*}
while $a =0 $,
\begin{eqnarray*}
 |\chi_j (\mu ) | &=&
         \frac {\frac 1{n^2 } \delta \pi +o(\frac 1{n^2 })}{
               (n+o(1))(|\theta|+o(1))} = \frac 1{n^3 } \frac {\delta
                \pi }{|\theta| } , \quad \mbox{ if $ \theta \ne a
          $, } \\  
 |\chi_j (\mu ) | &=&
         \frac {\frac 1{n^2 } \delta \pi +o(\frac 1{n^2 })}{
               (n+o(1))(\frac 1n \delta \pi +o(\frac 1n))} = \frac 1{n^2 }
          +o(\frac 1{n^2 }), \mbox{ if $ \theta = a $. }
\end{eqnarray*}
Summarize above results as : for $n$ sufficiently large,  
\begin{eqnarray*}
  || {\cal R}(\mu ) \frac {\sin \mu \pi }{\mu } ||_\infty &\sim & \mbox{}
        \frac 1{\delta n} |a| \pi , \quad 1cm \mbox{ if $ a \neq 0
$, } \\ 
|| {\cal R}(\mu ) \frac {\sin \mu \pi }{\mu } ||_\infty &\sim & \mbox{}
        O(\frac 1{n^2 } ), \quad \mbox{ if $ \alpha = 0 $. }
\end{eqnarray*}

This tells us that $|| {\cal R}(\mu ) \frac {\sin \mu \pi }{\mu } ||_\infty
$ has growth order at most $O(1/n)$ as $n$ is sufficiently large. With such
estimate oon the contours, we have the desired estimate on $|| U^* \Psi ^{-1}
 (\mu ) {\cal E} (\mu ) U || $, so we complete the proof {\bf Lemma \ref{Th2.2} 
}. 

%
%

\end{document}